\documentclass[onecolumn,notitlepage,12pt]{article}
\usepackage{amsfonts}
\usepackage{amsmath}
\usepackage{sectsty}
\usepackage{latexsym}
\usepackage{amsbsy}
\usepackage{amssymb}
\usepackage{verbatim}
\usepackage{amsthm}
\usepackage{makeidx}
\usepackage{graphicx}
\usepackage{layout}
\usepackage{pb-diagram}
\usepackage{amscd}
\usepackage{anysize}
\usepackage{tikz}
\usepackage{upgreek}
\usepackage[all,cmtip]{xy}
\usepackage{youngtab}
\usepackage{hyperref}
\usepackage{color}
\usepackage{mathrsfs}

\newtheorem{theorem}{Theorem}[section]

\newtheorem{claim}[theorem]{Claim}

\newtheorem{conjecture}[theorem]{Conjecture}
\newtheorem{corollary}[theorem]{Corollary}

\newtheorem{lemma}[theorem]{Lemma}

\newtheorem{proposition}[theorem]{Proposition}
\newtheorem{remark}[theorem]{Remark}

\def\ra{\rightarrow }
\def\M{\mathcal{M}}
\def\R{\mathcal{R}}
\def\Z{\mathbb{Z}}

\def\bthm{\begin{theorem}}
\def\ethm{\end{theorem}}
\def\bcl{\begin{claim}}
\def\ecl{\end{claim}}
\def\bcon{\begin{conjecture}}
\def\econ{\end{conjecture}}
\def\bpf{\begin{proof}}
\def\epf{\end{proof}}
\def\blm{\begin{lemma}}
\def\elm{\end{lemma}}
\def\bpp{\begin{proposition}}
\def\epp{\end{proposition}}
\def\brm{\begin{remark}}
\def\erm{\end{remark}}

\begin{document}
\title{On the moduli space of nonpositively curved metrics on a product with a torus}
\author{F. Thomas Farrell \and Yi Jiang}
\date{}
\maketitle

\begin{abstract}
We prove that if $N$ is a closed smooth manifold admitting a
Riemannian metric of negative sectional curvature and $T^k$ is the $k$-dimensional torus, then there is a
retraction from the moduli space of nonpositively curved
Riemannian metrics on $T^k\times N$ onto the moduli space of flat
metrics on $T^k$. As an application, we combine our retraction with Tuschmann and Wiemeler’s work on the rational topology of the moduli space of flat metrics on $T^k$ to deduce that the moduli space of
nonpositively curved metrics on $T^k\times N$ has nontrivial higher
rational homotopy and cohomology groups when $k\geq5$.
\end{abstract}

\section{Introduction}

Let $M$ be a closed smooth manifold and let $\mathcal R(M)$ denote
the space of Riemannian metrics on $M$, endowed with the
$C^\infty$-topology. The group $\operatorname{Diff}(M)$ of
self-diffeomorphisms of $M$ acts on $\mathcal R(M)$ by pulling back
Riemannian metrics. The quotient space
\[
\mathcal M(M):=\mathcal R(M)/\operatorname{Diff}(M)
\]
is called \emph{the moduli space of Riemannian metrics on $M$}.
Similarly, we define \emph{the moduli space of nonpositively curved metrics on $M$}
\[
\mathcal M_{\leq0}(M)
:=\mathcal R_{\leq0}(M)/\operatorname{Diff}(M),
\]
where $\mathcal R_{\leq0}(M)$ is the subspace of $\mathcal R(M)$
consisting of Riemannian metrics of nonpositive sectional curvature. We likewise define
$\mathcal M_{=0}(M):=\mathcal R_{=0}(M)/\operatorname{Diff}(M)$,
where $\mathcal R_{=0}(M)$ consists of all flat metrics on $M$.

Tuschmann and Wiemeler \cite[Theorems~1.1 and 5.4]{TW} obtained retraction results for
moduli spaces of Riemannian metrics of nonnegative Ricci curvature on products
with tori, and used them to detect nontrivial higher rational homotopy
and cohomology groups. In this paper, we prove an analogous retraction
result for moduli spaces of nonpositively curved metrics. We denote the
$k$-dimensional torus by $T^k$ throughout the paper.

\begin{theorem}\label{th-retraction}
Let $N$ be a closed smooth manifold admitting a Riemannian metric of
negative sectional curvature. Then there is a retraction
\[
r:\mathcal M_{\leq0}(T^k\times N)\longrightarrow
\mathcal M_{=0}(T^k).
\]
\end{theorem}

Although Theorem~\ref{th-retraction} is formally analogous to \cite[Theorem~5.4]{TW}, the two proofs rely on different geometric mechanisms. Tuschmann and Wiemeler apply the Cheeger--Gromoll splitting theorem \cite{CG} and construct their retraction by taking the flat base of a suitable Riemannian submersion; their continuity argument employs fiber shrinking and Gromov--Hausdorff convergence. In contrast, we combine the de~Rham decomposition of the universal cover with the classical theory of centers and Clifford translations \cite{LY,E82}. The resulting action of $\pi_1(T^k)\times1$ on the Euclidean factor yields the flat torus defining $r$. Continuity is proved by reducing smooth convergence of representative Riemannian metrics first to convergence of translation lengths and then, via polarization, to convergence of the associated lattice Gram matrices. The key point is that central deck transformations in nonpositive curvature are Clifford translations (Lemma~\ref{lm-Center-Clifford}), so their displacement functions are constant and equal to the lengths of the corresponding lattice vectors (Lemma~\ref{displace-lattice}). In particular, the metric on $T^k\times N$ need not be a Riemannian product.

By \cite[Proposition~5.5]{TW} and its correction \cite[Proposition~1]{TWC}, $\mathcal M_{=0}(T^k)$ is simply
connected and has nontrivial higher rational cohomology and
homotopy groups when $k\geq5$. We therefore obtain the following consequence.

\begin{corollary}\label{cor-main}
  Let $N$ be a closed smooth manifold admitting a Riemannian metric of negative sectional curvature and let $k\geq5$. Then $\M_{\leq0}(T^k\times N)$ has nontrivial higher rational cohomology groups and nontrivial higher rational homotopy groups.
\end{corollary}

After defining $r$ and verifying that it is well defined in Section~2, we establish its continuity in Section~3. The proof of Theorem~\ref{th-retraction} is completed in Section~4.

\section{Definition of the retraction \texorpdfstring{$r$}{r}}\label{se-define r}

Define a map $r:\M_{\leq0}(T^k\times N)\ra \M_{=0}(T^k)$ as follows. Let $h$ be a nonpositively curved Riemannian metric on $T^k\times N$. We construct an isometry class $[h']\in \M_{=0}(T^k)$ from $h$ and define $r([h]):=[h']$. Let $(\tilde{M},\tilde{h})$ be the universal Riemannian cover of $(T^k\times N,h)$. By the de~Rham decomposition \cite[p.12]{E}, $(\tilde{M},\tilde{h})$ can be written as a Riemannian product
\[
(\tilde{M},\tilde{h})=E\times Q
\]
where $E$ is the Euclidean de~Rham factor and $Q$ is the product of all non-Euclidean de~Rham factors. Since the group $\Gamma\cong\pi_1(T^k\times N)$ of covering transformations is a discrete subgroup of $I(\tilde{M},\tilde{h})=I(E)\times I(Q)$ \cite[p.13, p.368 Remark]{E}, it follows from \cite[Proposition~1.9.4, Theorem~7.2.1]{E} that the subgroup $C(\Gamma)$ of Clifford translations in $\Gamma$  is embedded as a discrete subgroup of $I(E)$ by the projection \[
C(\Gamma)\subset I(\tilde{M},\tilde{h})\ra I(E)
\]
and the quotient space $E/C(\Gamma)$ is compact. 
Furthermore, $C(\Gamma)$ is equal to the center $Z(\Gamma)$ of $\Gamma$ by Lemma~\ref{lm-Center-Clifford} below. Since $\pi_1(N)$ is centerless\footnote{This is implicitly proved in \cite[Section~9]{EO}. See \cite[Remark~7]{BFJ} for an outline of the argument.}, it follows that $Z(\Gamma)\cong\pi_1(T^k)\cong\Z^k$ is torsion-free and hence $Z(\Gamma)$ acts freely on $E$. This implies that $E/Z(\Gamma)$ is a compact flat manifold which is diffeomorphic to $T^k$ by the Bieberbach theorem \cite[p.65, Theorem~5.4]{C}. Thus $E/Z(\Gamma)$ represents an isometry class $[h']\in \M_{=0}(T^k)$.

\blm\label{lm-Center-Clifford} Let $M=T^k\times N$, where $N$ is a closed smooth manifold admitting a Riemannian metric of negative sectional curvature, and let $g$ be a Riemannian metric on $M$ with nonpositive sectional curvature. Denote $\pi_1(M)$ by $\Gamma$ and regard it as a subgroup of the isometry group of the universal cover of $(M,g)$. Then $Z(\Gamma)=C(\Gamma)$, where $C(\Gamma)$ denotes the subgroup of Clifford translations in $\Gamma$.
\elm
\bpf
Since  $C(\Gamma)$ equals the unique maximal abelian normal subgroup of $\Gamma$ (cf.~\cite[p.210 and Section~1.4]{Ep}),
\begin{comment}cf.~\cite[p.210, Line -4]{Ep} and note that $Gamma$ satisfies the duality condition since $M=\tilde{M}/\Gamma$ is compact (cf.~\cite[1.2]{Ep})
In fact \cite[p.214 the last two sentences in section 1.4]{Ep} says that this can be deduced directly from \cite[Theorem~2.4]{CE} which is as follows:
      \begin{quote}
        Let $H$ be a Hadamard  manifold and let $G\subset I(H)$ be a subgroup with nontrivial center $A$ such that the normalizer $N(G)$ of $G$ in $I(H)$ satisfy the duality condition. Then $A$ consists of Clifford translations.
      \end{quote}
      Now let $H$ be the universal cover of $M$ and let $G$ be the maximal abelian normal subgroup of $\Gamma=\pi_1(M)$. Then $G=A$ and $\Gamma\subset N(G)$. Since $H/\Gamma$ is compact, $\Gamma$ satisfies the duality condition.\cite[p.213, Line 9-11]{Ep}. Since $\Gamma\subset N(G)$, the definition of duality condition implies that $N(G)$ satisfies the duality condition. Hence we can apply \cite[Theorem~2.4]{CE} here to get $A=G$ is contained in the Clifford subgroup of $\Gamma$.
\end{comment}
it suffices to show that $Z(\Gamma)$ is the maximal abelian normal subgroup of $\Gamma$. Let $H$ be any abelian normal subgroup of $\Gamma$. Then the image $q(H)$ of the composition \[
q:H\subset \Gamma\cong\pi_1(T^k)\times \pi_1(N)\ra \pi_1(N)
\]
            is an abelian normal subgroup of $\pi_1(N)$. By Preissmann's theorem, $q(H)$ is either trivial or infinite cyclic. Suppose that $q(H)$ is infinite cyclic. Regard $\pi_1(N)$ as the group of covering transformations on the universal cover $\tilde{N}$ of $N$ and let $C(\pi_1(N))$ denote the subgroup of Clifford translations in $\pi_1(N)\subset I(\tilde{N})$. Since $q(H)$ is an abelian normal subgroup of $\pi_1(N)$, the Clifford subgroup $C(\pi_1(N))$ contains $q(H)$ and hence has rank at least $1$. This contradicts the fact that any Clifford translation on a negatively curved simply connected manifold is trivial (cf.~\cite[p.32, Proposition~1.9.4]{E}). Thus $q(H)$ is trivial and hence $H\subset \pi_1(T^k)\times 1=Z(\Gamma)$. Consequently $Z(\Gamma)$ is the maximal abelian normal subgroup of $\Gamma$ and hence $C(\Gamma)=Z(\Gamma)$.
\epf

\begin{lemma}
  $r:\M_{\leq0}(T^k\times N)\ra \M_{=0}(T^k)$ is a well-defined map.
\end{lemma}
\bpf
Let $h_1,h_2$ be two nonpositively curved Riemannian metrics on $T^k\times N$ and let $f: (T^k\times N, h_1)\ra (T^k\times N, h_2)$ be an isometry. For $i=1,2$, the universal cover $(\tilde{M},\tilde{h}_i)$ can be written as a Riemannian product 
\[
(\tilde{M},\tilde{h}_i)=E_i\times Q_i
\]
where $E_i$ is a Euclidean space and $Q_i$ has no Euclidean factor. The isometry $f$ lifts to an isometry $\tilde{f}:(\tilde{M},\tilde{h}_1)\ra (\tilde{M},\tilde{h}_2)$ which splits as $\tilde{f}=(f_E,f_Q)$ where both $f_E:E_1\ra E_2$ and $f_Q:Q_1\ra Q_2$ are isometries\footnote{By the uniqueness of the distributions in the de~Rham decomposition \cite[p.185, Theorem~5.4(4)]{KN}, one has
\begin{gather*}
  d\tilde{f}(TE_1\times Q_1)\subset TE_2\times Q_2\\
 d\tilde{f}(E_1\times TQ_1)\subset E_2\times TQ_2
\end{gather*}
The criterion for an isometry to split \cite[p.13]{E} therefore implies that $\tilde{f}$ splits.}.
Since $\tilde{f}Z(\Gamma)\tilde{f}^{-1}\subset Z(\Gamma)$\footnote{Let $p:\tilde{M}\rightarrow M$ be the covering map and $p(\tilde{x_0})=x_0$. Let $\alpha_{\tilde{x_0}}:\pi_1(M,x_0)\rightarrow\Gamma$ be the group isomorphism defined as follows. For each loop $\gamma$ based at $x_0$, let $\tilde{\gamma}$ be the lift with $\tilde{\gamma}(0)=\tilde{x_0}$. Now define $\alpha_{\tilde{x_0}}$ to be the unique covering transformation of $\tilde{M}$ which maps $\tilde{x_0}$ to $\tilde{\gamma}(1)$. The formula 
\[
\tilde{f}\alpha^{\tilde{x_0}}([\gamma])\tilde{f}^{-1}=\alpha^{\tilde{f}(\tilde{x_0})}(f_*[\gamma])
\]
implies that $\tilde{f}Z(\Gamma)\tilde{f}^{-1}\subset Z(\Gamma)$.}, it follows that $f_EZ(\Gamma)f_E^{-1}\subset Z(\Gamma)$ and this shows that $f_E:E_1\rightarrow E_2$ induces an isometry $E_1/Z(\Gamma)\cong E_2/Z(\Gamma)$. Consequently, $r$ is well-defined.
\epf

\section{Continuity of the map \texorpdfstring{$r$}{r}}
In this section, we prove the continuity of $r$ in Lemma~\ref{continuous}. We begin with a well-known fact.
\begin{lemma}\label{metrizability of moduli space}
Let $M$ be a smooth closed manifold. Then there is a $\operatorname{Diff}(M)$-invariant metric $d$ on $\R(M)$ (i.e., $d(f^*h,f^*g)=d(g,h)$ for all $f\in\operatorname{Diff}(M)$ and $g,h\in\R(M)$) that induces the $C^\infty$-topology on $\R(M)$ and for which the function $\tilde{d}:\M(M)\times\M(M)\rightarrow\mathbb{R}$, defined by $\tilde{d}([g],[h])=\inf_{f\in\operatorname{Diff}(M)}d(f^*g,h)$, is a metric on $\M(M)$ and generates the quotient $C^{\infty}$-topology on $\M(M)$.
\end{lemma}

See \cite[Lemma~2.1]{BGA22} and \cite[Proposition~4.1]{FO10} for a proof. Both arguments are based on Ebin's slice theorem \cite{Ebin}. The proof of \cite[Proposition~4.1]{FO10} can be modified to give a proof of Lemma~\ref{metrizability of moduli space} by replacing $\operatorname{Diff}_0(M)$ with $\operatorname{Diff}(M)$.

\begin{corollary}\label{cor-convergence of seq}
  Suppose $[h_0],[h_m]\in\M(M)$. Then $[h_m]$ converges to $[h_0]$ if and only if, for each positive integer $m$, there is $\phi_m\in\operatorname{Diff}(M)$ such that $\phi_m^*h_m$ converges to $h_0$.
\end{corollary}
\begin{proof}
  The ``if'' part follows directly from the continuity of the quotient map $\R(M)\rightarrow\M(M)$. For the ``only if'' part, the definition $\tilde{d}([h_m],[h_0])=\inf_{f\in\operatorname{Diff}(M)}d(f^*h_m,h_0)$ implies that, for each positive integer $m$, there is $\phi_m\in\operatorname{Diff}(M)$ such that
  \[
0\leq d(\phi_m^*h_m,h_0)<\tilde{d}([h_m],[h_0])+\frac{1}{m}.
\]
  As $m\to\infty$, we have $d(\phi_m^*h_m,h_0)\to0$; that is, $\phi_m^*h_m$ converges to $h_0$.
\end{proof}
Let $h$ be a nonpositively curved Riemannian metric on $M=T^k\times N$, and let $(\tilde{M},\tilde{h})$ be the universal Riemannian cover of $(M,h)$. Recall that $(\tilde{M},\tilde{h})$ can be written as a Riemannian product
\[
(\tilde{M},\tilde{h})=E\times Q,
\]
where $E$ is the Euclidean de~Rham factor and $Q$ is the product of all non-Euclidean de~Rham factors. For $\gamma\in Z(\Gamma)$, let $\gamma_E$ denote its image under the composition $Z(\Gamma)\subset\Gamma\subset I(E\times Q)\ra I(E)$. Then $\gamma_E$ is an ordinary translation of $E$, and the image of $Z(\Gamma)$ under the composition $Z(\Gamma)\subset\Gamma\subset I(E\times Q)\ra I(Q)$ is trivial by \cite[Proposition~1.9.4]{E}. Let $v:Z(\Gamma)\ra E$ be defined by $v(\gamma)=\gamma_E(e)-e$, where $e\in E$; this definition is independent of the choice of $e\in E$.

\blm\label{displace-lattice} The map $v:Z(\Gamma)\ra E$ is a homomorphism. Moreover, for each $\gamma\in Z(\Gamma)$, the displacement function $d_{\tilde{h}}(x,\gamma x)$ is constant and equal to the length $l(\gamma):=\|v(\gamma)\|$.
\elm
\bpf
For any $\gamma_1,\gamma_2\in Z(\Gamma)$ and any $(e,q)\in E\times Q$, \cite[Proposition~1.9.4]{E} implies that $(\gamma_1\gamma_2)(e,q)=(e+v(\gamma_1\gamma_2),q)$. On the other hand,
\[
\gamma_1(\gamma_2(e,q))=\gamma_1(e+v(\gamma_2),q)=(e+v(\gamma_2)+v(\gamma_1),q).
\]
Thus $v(\gamma_1\gamma_2)=v(\gamma_1)+v(\gamma_2)$.

Moreover, for any $x=(e,q)\in E\times Q$, we have
\[
d_{\tilde{h}}(x,\gamma x)=d_{\tilde{h}}((e,q),(e+v(\gamma),q))=\|v(\gamma)\|.
\]
\epf
We call $v$ the \emph{lattice homomorphism} for the Riemannian metric $h$. The image of $v$ is the lattice of the flat torus $E/Z(\Gamma)$ representing $r([h])$.

To prove the continuity of $r$, let us recall the following basic fact on $\M_{=0}(T^k)$.
$\M_{=0}(T^k)$ with the $C^0$-topology is homeomorphic to $O(k)\backslash\operatorname{GL}(k,\mathbb{R})/\operatorname{GL}(k,\Z)$ with the biquotient topology. The homeomorphism sends the metric of a flat torus $E/\Lambda$ to its representation matrix in standard affine coordinates. Its inverse sends $A\in\operatorname{GL}(k,\mathbb{R})$ to the flat torus $T_A=\mathbb{R}^k/\Z^k$ equipped with the pullback metric under the map $A:\mathbb{R}^k/\Z^k\ra\mathbb{R}^k/A\Z^k$; in canonical affine coordinates, the metric matrix is $A^TA$ (cf.~\cite[Corollary~of Theorem~1]{W73} and \cite[Section~1.3]{Garcia19}). Since the representation matrix is constant and its higher derivatives vanish, the $C^0$- and $C^\infty$-topologies coincide on $\M_{=0}(T^k)$.

\blm\label{continuous} Let $N$ be a closed smooth manifold admitting a Riemannian metric of negative sectional curvature. Then the map
\[
r:\M_{\leq0}(T^k\times N)\ra \M_{=0}(T^k)
\]
defined in Section~\ref{se-define r} is continuous.
\elm
\begin{comment}
Recall the definition of $r$ as follows. Let $h$ be a nonpositively curved Riemannian metric on $T^k\times N$. Then the universal cover $\widetilde{(T^k\times N,h)}$ of $(T^k\times N,h)$ is decomposed as a Riemannian product:
\[
\widetilde{(T^k\times N,h)}=E\times Q
\]
where $E$ is the Euclidean de~Rham factor and $Q$ is the product of all non-Euclidean de~Rham factors. Let $Z$ denote the center of $\pi_1(T^k\times N)$ which is identified with a subgroup of the covering transformation group and hence is a discrete subgroup of the isometry group $I(E\times Q)=I(E)\times I(Q)$. The composition of homomorphisms 
\[
\xymatrix{Z^k\cong Z\ar[r]^-{\mbox{incl.}}&I(E\times Q)\ar[r]^-{\mbox{proj.}}&I(E)}
\]
is monic and $E/Z$ is isomorphic to a flat torus $(T^k,g)$. Define $r([h])$ to be [g].
\end{comment}
\bpf
Since $\M_{\leq0}(T^k\times N)$ is first countable by Lemma~\ref{metrizability of moduli space}, it suffices to prove the sequential continuity of $r$. Let $\{[h_m]\}_m$ be a sequence converging to $[h_0]$ in $\M_{\leq0}(T^k\times N)$. We must show that $r([h_m])\to r([h_0])$ as $m\to\infty$.

By Corollary~\ref{cor-convergence of seq}, $[h_m]\to[h_0]$ if and only if there is a sequence of diffeomorphisms $\phi_m\in\operatorname{Diff}(T^k\times N)$ such that $g_m:=\phi_m^*h_m\to h_0$ as $m\to\infty$.

\begin{comment}
Let $(\tilde{M},\tilde{h}_m)$ be the universal cover of $(M,h_m)$ and let $\tilde{\phi}_m\in\operatorname{Diff}(\tilde{M})$ be the lifting diffeomorphism of $\phi_m$.
Consider the following commutative diagram:
\[
\xymatrix{(\tilde{M},\tilde{\phi}_m^*h_m)\ar[r]^-{\tilde{\phi}_m}_-{\mbox{isom.}}\ar[d]&(\tilde{M},\tilde{h}_m)\ar[d]\\
(M,\phi_m^*h_m)\ar[r]^-{\phi_m}_-{\mbox{isom.}}&(M,h_m)}
\]
It follows that
\[
\xymatrix{\tilde{g}_m:=\tilde{\phi}_m^*\tilde{h}_m\ar[r]_-{m\ra\infty}&\tilde{h}_0}.
\]
(under uniformly convergence on compact sets, Please check the expression?)
\end{comment}

Let $(\tilde{M},\tilde{g}_m)$ be the universal cover of $(M,g_m)$. For any $\gamma\in Z(\Gamma)=\pi_1(T^k)\times 1\subset\Gamma\subset I(\tilde{M},\tilde{g}_m)$, Lemma~\ref{lm-Center-Clifford} implies that $\gamma$ is a Clifford translation for the Riemannian metric $\tilde{g}_m$. Hence the displacement function $d_{\tilde{g}_m}(x,\gamma x)$ is constant and, by Lemma~\ref{displace-lattice}, equals the length $l_m(\gamma)$ of the lattice vector $v_m(\gamma)$ in the flat torus $r([h_m])$. Since $g_m\to h_0$ as $m\to\infty$, it follows from Lemma~\ref{lm-displace inequality} below that $l_m(\gamma)$ converges to $l_0(\gamma)$ for every $\gamma\in Z(\Gamma)$.

Let
\[
Q_m(\gamma):=l_m(\gamma)^2=\langle v_m(\gamma),v_m(\gamma)\rangle_m
\]
where $\langle\cdot,\cdot\rangle_m$ denotes the inner product on the Euclidean de~Rham factor $E_m$ of $(\tilde{M},\tilde{g}_m)$. Since $l_m(\gamma)$ converges to $l_0(\gamma)$, it follows that $Q_m(\gamma)$ converges to $Q_0(\gamma)$ for every $\gamma\in Z(\Gamma)$.

Let $e_1,\ldots,e_k$ be a standard basis of $Z(\Gamma)=\pi_1(T^k)\times 1\cong\Z^k$, and consider the Gram matrix $G_m:=(\langle v_m(e_i),v_m(e_j)\rangle_m)_{ij}$. Since $v_m$ is a homomorphism by Lemma~\ref{displace-lattice}, we have
\[
(G_m)_{ij}=\frac{1}{2}[Q_m(e_i+e_j)-Q_m(e_i)-Q_m(e_j)]
\]
for every nonnegative integer $m$. Hence $(G_m)_{ij}$ converges to $(G_0)_{ij}$ for all $i,j$, and therefore $G_m$ converges to $G_0$ in the $C^0$-topology. The Gram matrix represents the metric of a flat torus in standard affine coordinates. It is constant, and all its higher derivatives vanish; hence $[G_m]$ converges to $[G_0]$ in the $C^\infty$-topology on $\M_{=0}(T^k)$. That is, $r([h_m])$ converges to $r([h_0])$ in the $C^\infty$-topology on $\M_{=0}(T^k)$.
\epf

To prove Lemma~\ref{lm-displace inequality}, we first establish the following lemma.
\begin{lemma}\label{lm-linear algebra}
  Let $T$ be a symmetric bilinear form on a Euclidean space $V$. Then
  \[
|T(v,v)|\leq\|T\|\|v\|^2
\]
  for every $v\in V$, where $\|T\|^2=\sum_{i,j}T(e_i,e_j)^2$ for an orthonormal basis $\{e_i\}$ of $V$; this norm is independent of the choice of basis.
\end{lemma}
\begin{comment}
\brm Verify that $\|T\|^2$ does not depend on the choice of $\{e_i\}$: Let $\{e_i\}$ and $\{e'_i\}$ be two orthonormal bases and let $A$ be the transition matrix, i.e.,
\[
(e_1,\cdots,e_n)^T=A(e'_1,\cdots,e'_n)^T
\]
and hence $A$ is an orthogonal matrix.
Let $B:=(T(e_i,e_j))_{ij}$ and $B':=(T(e'_i,e'_j))_{ij}$. Then $B=AB'A^T$ and hence 
\begin{eqnarray*}
 \sum_{i,j}T(e_i,e_j)^2 &=& tr(B^TB) \\
  &=& tr(A(B')^TA^TAB'A^T) \\
  &=& tr(A(B')^TB'A^{-1})  \\
  &=& tr((B')^TB')=\sum_{i,j}T(e'_i,e'_j)^2 
\end{eqnarray*}
\erm
\end{comment}
\bpf
Let $\hat{T}:V\ra V$ be the self-adjoint operator associated with $T$, i.e., $T(v,w)=\langle\hat{T}v,w\rangle$ for all $v,w\in V$. Then there is an orthonormal basis $\{e_i\}$ of eigenvectors of $\hat{T}$. Let $\lambda_i$ be the eigenvalue associated with $e_i$. Then
\[
\|T\|^2=\sum_{i,j}T(e_i,e_j)^2=\sum_{i,j}\langle\hat{T}e_i,e_j\rangle^2=\sum_{i,j}\langle\lambda_ie_i,e_j\rangle^2=\sum_{i}\lambda_i^2.
\]

For any $v\in V$, let $v=\sum_{i}v_ie_i$. Then \[
|T(v,v)|=|\langle\hat{T}v,v\rangle|=|\langle\sum_{i}\lambda_iv_ie_i,\sum_jv_je_j\rangle|=|\sum_i\lambda_iv_i^2|
\]
and hence \[
|T(v,v)|\leq\sum_i|\lambda_i|v_i^2\leq\max_i|\lambda_i|\sum_iv_i^2=\max_i|\lambda_i|\|v\|^2\leq\sqrt{\sum_i\lambda_i^2}\|v\|^2=\|T\|\|v\|^2.
\]
\epf

Let $N$ be a closed smooth manifold admitting a Riemannian metric of negative sectional curvature, and let $\{g_m\}_m$ and $h_0$ be Riemannian metrics of nonpositive sectional curvature on $M=T^k\times N$. Denote the lattice homomorphism for $g_m$ by $v_m$ and that for $h_0$ by $v_0$. Let $\Gamma=\pi_1(M)$.

\begin{lemma}\label{lm-displace inequality}Suppose that $g_m$ converges to $h_0$ in the $C^\infty$-topology as $m\to\infty$. Then the length $l_m(\gamma)$ of $v_m(\gamma)$ converges to the length $l_0(\gamma)$ of $v_0(\gamma)$ for every $\gamma\in Z(\Gamma)$.
\end{lemma}
\bpf We first show that, for every $\varepsilon\in(0,1)$, there is a positive integer $m_0$ such that
  \begin{equation}\label{ine-displace}
    \sqrt{1-\varepsilon}d_0(x,\gamma x)\leq d_m(x,\gamma x)\leq\sqrt{1+\varepsilon}d_0(x,\gamma x)
  \end{equation}
for all $m>m_0$ and all $x\in\tilde{M}$, where $d_m$ denotes $d_{\tilde{g}_m}$ and $d_0$ denotes $d_{\tilde{h}_0}$.
  Let $S_m:=g_m-h_0$. Since $g_m$ converges to $h_0$, we have $\sup_{x\in M}\|S_m\|_{h_0}(x)\to0$, where
  \[
\|S_m\|_{h_0}(x):=\left(\sum_{i,j}S_m(e_i,e_j)^2\right)^{1/2}
\]
and $\{e_i\}$ is an $h_0$-orthonormal basis at $x$. Thus, for every $\varepsilon>0$, there is a positive integer $m_0$ such that $\sup_{x\in M}\|S_m\|_{h_0}(x)\leq\varepsilon$ for all $m>m_0$. By Lemma~\ref{lm-linear algebra}, we have
\[
|S_m(v,v)|\leq\|S_m\|_{h_0}(x)h_0(v,v)
\]
for all $x\in M$ and $v\in T_xM$. Hence, for every $\varepsilon>0$, there is a positive integer $m_0$ such that
\[
|S_m(v,v)|\leq\varepsilon h_0(v,v)
\]
for all $m>m_0$; that is, $|g_m(v,v)-h_0(v,v)|\leq\varepsilon h_0(v,v)$, and therefore
\[
(1-\varepsilon)h_0(v,v)\leq g_m(v,v)\leq(1+\varepsilon)h_0(v,v)
\]
for all $m>m_0$.
Pulling back these Riemannian metrics to $\tilde{M}$ gives 
\begin{equation}\label{ine-metric}
(1-\varepsilon)\tilde{h}_0(\tilde{v},\tilde{v})\leq \tilde{g}_m(\tilde{v},\tilde{v})\leq(1+\varepsilon)\tilde{h}_0(\tilde{v},\tilde{v})
\end{equation}
for every tangent vector $\tilde{v}$ to $\tilde{M}$.
\begin{comment}
The detail is as follows:
For any $\tilde{x}\in\tilde{M}$ and $\tilde{v}\in T_{\tilde{x}}\tilde{M}$ we have
\[
\tilde{g}_m(\tilde{v},\tilde{v})=g_m(d\pi(\tilde{v}),d\pi(\tilde{v}))
\]
and the same for $\tilde{h}_0$. Then
\[
(1-\varepsilon)h_0(d\pi(\tilde{v}),d\pi(\tilde{v}))\leq g_m(d\pi(\tilde{v}),d\pi(\tilde{v}))\leq(1+\varepsilon)h_0(d\pi(\tilde{v}),d\pi(\tilde{v}))
\]
and this implies the inequality \eqref{ine-metric}.
\end{comment}
For any Riemannian metric $g$ on $\tilde{M}$, let $d_g$ denote the associated distance function. For any $x,y\in\tilde{M}$, let $L(x,y)$ be the set of all smooth curves connecting $x$ and $y$. Then
\[
d_g(x,y)=\inf_{l\in L(x,y)}\int_{0}^{1}\sqrt{g(l'(t),l'(t))}\,dt.
\]
Thus inequality \eqref{ine-metric} implies inequality \eqref{ine-displace}.

For each fixed $\gamma\in Z(\Gamma)$, inequality \eqref{ine-displace} and Lemma~\ref{displace-lattice} imply that, for every $\varepsilon\in(0,1)$, there is a positive integer $m_0$ such that
\[
\sqrt{1-\varepsilon}l_0(\gamma)\leq l_m(\gamma)\leq\sqrt{1+\varepsilon}l_0(\gamma)
\]
for all $m>m_0$. Letting $\varepsilon\to0$ yields $\lim_{m\to\infty}l_m(\gamma)=l_0(\gamma)$.
\epf

\section{Proof of Theorem~\ref{th-retraction}}
\bpf[Proof of Theorem~\ref{th-retraction}]
By Lemma~\ref{continuous}, $r$ is continuous. Moreover, $r$ has a continuous section defined by taking the Riemannian product with $(N,g)$, where $g$ is a fixed negatively curved metric on $N$. This completes the proof.
\epf

\section*{Acknowledgements}

The authors thank Wilderich Tuschmann for very helpful discussions on the proofs of Lemmas \ref{metrizability of moduli space} and \ref{continuous}.
Yi Jiang is partially supported by NSFC12371070.

\section*{AI Declaration}

Yi Jiang used Kimi K3 (Moonshot AI) to review the mathematical arguments in an earlier draft. Kimi suggested a simplification of the proof of Lemma~\ref{continuous}, allowing an assumption in the earlier version to be removed. Both authors carefully verified the suggestion and wrote the final proof themselves. Kimi was also used for English-language polishing and for citation and LaTeX formatting; the authors reviewed and approved all resulting changes. The authors take full responsibility for the content of the manuscript.

{\footnotesize
\bibliographystyle{alpha}
\bibliography{bib}
}

\bigskip
\noindent
F. Thomas Farrell, 1200 Elm ST Unit 529, Manchester, NH 03101\\
\texttt{ftfarrell@gmail.com}

\bigskip
\noindent
Yi Jiang, Academy for Multidisciplinary Studies, Beijing National Center
for Applied Mathematics, Capital Normal University, Beijing, 100048, China\\
\texttt{jiangyi@cnu.edu.cn}

\end{document}